\documentclass[12pt,onecolumn,twoside]{IEEEtran}

\usepackage{amsmath,amssymb,euscript,psfrag,latexsym,color,graphicx}

\begin{document}
\newtheorem{thm}{Theorem}
\newtheorem{cor}[thm]{Corollary}
\newtheorem{lemma}[thm]{Lemma}
\newtheorem{prop}[thm]{Proposition}
\newtheorem{problem}[thm]{Problem}
\newtheorem{remark}[thm]{Remark}
\newtheorem{defn}[thm]{Definition}
\newtheorem{ex}[thm]{Example}
\newcommand{\ignore}[1]{}
\newcommand{\mR}{{\mathbb R}}
\newcommand{\mN}{{\mathbb N}}
\newcommand{\bR}{{\bf R}}
\newcommand{\rR}{{\rm R}}
\newcommand{\mZ}{{\mathbb Z}}
\newcommand{\mC}{{\mathbb C}}
\newcommand{\fR}{{\mathfrak R}}
\newcommand{\fK}{{\mathfrak K}}
\newcommand{\fG}{{\mathfrak G}}
\newcommand{\fT}{{\mathfrak T}}
\newcommand{\fC}{{\mathfrak C}}
\newcommand{\fF}{{\mathfrak F}}
\newcommand{\cF}{{\mathcal F}}
\newcommand{\fM}{{\mathfrak M}}
\newcommand{\mD}{{\mathbb{D}}}
\newcommand{\cU}{{\mathcal{U}}}
\newcommand{\cL}{{\mathcal{L}}}
\newcommand{\cD}{{\mathcal{D}}}
\newcommand{\bJ}{{\mathbb{J}}}
\newcommand{\bI}{{\mathbb{I}}}
\newcommand{\s}{{\rm s}}
\newcommand{\ar}{{\rm a}}
\newcommand{\intpi}{{\frac{1}{2\pi}\int_{-\pi}^\pi}}
\newcommand{\cE}{{\mathcal E}}
\newcommand{\me}{{\rm _{ME}}}
\newcommand{\mr}{{\rm _{sm}}}

\title{Distances between power spectral densities$^*$}
\author{Tryphon T. Georgiou$^\dag$\thanks{$^\dag$Department of Electrical and Computer Engineering, University of Minnesota, Minneapolis, MN 55455; {\tt tryphon@ece.umn.edu}. $^*$This research has been supported by the NSF and the AFOSR.}}% , {\em IEEE Fellow}}
\date{}
\markboth{July 2006}
{Georgiou: power spectral densities \& distance measures}
\maketitle
\begin{abstract} We present several natural notions of distance between spectral density functions of (discrete-time) random processes. They are motivated by certain filtering problems. First we quantify the degradation of performance of a predictor which is designed for a particular spectral density function and then it is used to predict the values of a random process having a different spectral density. The logarithm of the ratio between the variance of the error, over the corresponding minimal (optimal) variance, produces a measure of distance between the two power spectra with several desirable properties.  Analogous quantities based on smoothing problems produce alternative distances and suggest a class of measures based on fractions of generalized means of ratios of power spectral densities. These distance measures endow the manifold of spectral density functions with a (pseudo) Riemannian metric. We pursue one of the possible options for a distance measure, characterize the relevant geodesics, and compute corresponding distances. \end{abstract}


\section{Introduction}
\PARstart{D}{ue} the centrality of spectral analysis in a wide range of scientific disciplines, there has been a variety of viewpoints regarding how to quantify distances between spectral density functions. Besides the obvious ones which are based on norms, inherited by ambient function-spaces $L_2,L_1$, etc., there has been a plethora of alternatives which attempt to acknowledge the structure of power spectral density functions as a positive cone. The most most well known are the Kullback-Leibler divergence which originates in hypothesis testing and in Bayes' estimation, the Itakura-Saito distance which originates in speech analysis ---both belonging to Bregman class (\cite{Gray,Bregman,GL}), the Bhattacharyya distance \cite{bhattacharyya}, and the Ali-Silvey class of divergences \cite{alisilvey}. Their origin can be traced either to a probabilistic rationale (as in the case of the Kullback-Leibler divergence) or, to some ad-hoc mathematical construct designed to seek distance measures with certain properties (as in the case of Bregman and Ali-Silvey classes). The purpose of this work is to introduce certain new notions of distance which are rooted in filtering theory and provide intrinsic distance measures between any two {\em power spectral density} functions. 

Our starting point is a prediction problem. We select an optimal predictive filter for an underlying random process based on the assumption that the process has a given power spectral density $f_1(\theta)$. We then evaluate the performance of such a filter against a second power spectral density $f_2(\theta)$ ---which may be thought of as the spectral density function of the ``actual'' random process. The relative degradation of performance (i.e., variance of the prediction error) quantifies a mismatch between the two functions. Interestingly, it turns out to be equal to the ratio of the arithmetic over the geometric mean of the fraction of the two  power spectra. The logarithm of the relative degradation serves a distance measure.  

Infinitesimal analysis suggests a pseudo-Riemannian metric on the manifold of power spectral density functions. The presence of such a metric suggests that geodesic distances may be used to quantify divergence between power spectra. Indeed, a characterization of geodesics is provided, and certain logarithmic intervals are shown to satisfy the condition. The length of such intervals connecting two power spectral densities  provides yet another notion of distance between the two.

An identical approach based on the degradation of performance of smoothing filters leads to other expressions which, equally well, quantify divergence between power spectral densities. Two observations appear to be universal. First that the mismatch between the ``shapes'' of spectral density functions is what turns out to be important. This is quantified by how far the ratio $f_1/f_2$ of the two spectral densities is from being constant across frequencies. The ratio of spectral density functions is reminiscent of the likelihood ratio in probability theory.

The second observation is that all of the distance measures that we encountered, in essence, they compare different means (i.e., arithmetic, geometric, and harmonic, possibly, weighted) of the ratio of the two spectral density functions or of their logarithms.
It is quite standard, that e.g., the argithmetic and the geometric means coincide only when the ratio is constant and have a gap otherwise.  The same applies to a wider family of generalized means. Thus, this observation suggests a much larger class of possible alternatives: quantify the divergence between (the ``shape'' of) two density functions using the gap between two generalized means of their ratio, or by the slackness of Jensen-type of inequalities involving this ratio. The underlying mathematical construct appears quite distinct from those utilized in defining the Bregman and the Ali-Silvey classes of distance measures. Furthermore, the mathematical construct is deeply rooted in prediction theory and, at least in certain cases, can be motivated as quantifying degradation of performance as we explained earlier.

\section{Preliminaries on least-variance prediction and smoothing}\label{section:preliminaries}

Consider a scalar zero-mean stationary random process
$\{u_k,\;k\in\mZ\}$ and denote by
$R_0,\,R_1,\,R_2,\ldots$ its sequence of autocorrelation samples and by
$d\mu(\theta)$ its power spectrum. Thus, $R_k:=\cE\{u_\ell u_{\ell+k}^*\}=R_{-k}^*$ and
\[
R_k=\intpi e^{jk\theta}d\mu(\theta),\mbox{ for }k\in\mZ,
\]
while $\cE$ denotes expectation and ``$^*$'' denotes complex conjugation.
We are interested in quadratic optimization problems with respect to the usual inner product
\begin{eqnarray}\nonumber
\langle \sum_k a_ku_k,\sum_\ell b_\ell u_\ell\rangle &:=&\cE\{ (\sum_k a_ku_k)(\sum_\ell b_\ell u_\ell)^*\}\\
&=&\sum_{k,\ell} a_kR_{k-\ell}b_\ell^*. \label{Tproduct}
\end{eqnarray}

The closure of ${\rm span}\{u_k\,:\,k\in\mZ\}$, which we denote by $\cU$, can be identified with the space $L_{2, d\mu}[-\pi,\pi)$ of functions which are square integrable with respect to $d\mu(\theta)$ with inner product
\[
\langle a,b\rangle_{d\mu}:=
\frac{1}{2\pi}\int_{-\pi}^\pi a(\theta)(b(\theta))^*d\mu(\theta)
%\langle \sum_k a_kz^k,\sum_\ell b_\ell z^\ell\rangle_{d\mu}:=
%\int_{-\pi}^\pi(\sum_k a_ke^{jk\theta})(\sum_\ell b_\ell e^{j\ell\theta})^*\frac{d\mu(\theta)}{2\pi}.
\]
where $a(\theta)=\sum_k a_ke^{jk\theta}$ and $b(\theta)=\sum_\ell b_\ell e^{j\ell\theta}$.
%It can also be identified with square-integrable series where their respective inner product and norm are defined via (\ref{Tproduct}).
Further, the correspondence
\begin{eqnarray*}
\cU\to L_{2,d\mu}&:& \sum_k a_ku_k \mapsto \sum_k a_ke^{jk\theta}
%\cU\to \ell_{2,\bR} &:& \sum_k a_ku_k \mapsto\left(\ldots,\,a_{-1},\,a_0,\,a_1,\ldots\right)
\end{eqnarray*}
is a Hilbert space isomorphism (see \cite{masani}).
Thus, least-variance approximation problems can be equivalently expressed in $L_{2,d\mu}$.
In particular, the variance of the {\em one-step-ahead prediction} error
$u_0-\hat{u}_{0|{\rm past}}$
for the {\em predictor}
\[\hat{u}_{0|{\rm past}}=\sum_{k>0}\alpha_ku_{-k}\] is
\begin{equation}\label{infalphas}
\cE\{|u_0 -\hat{u}_{0|{\rm past}}|^2\}= \|1-\sum_{k>0}\alpha_ke^{jk\theta}\|_{d\mu}^2.
\end{equation}
Similarly, the variance of the error
of the {\em smoothing} filter
\begin{equation}\label{uboth}
\hat{u}_{0|{\rm past \;\&\; future}}:=\sum_{k\neq 0} \beta_k u_{-k}
\end{equation}
is simply
\begin{equation}\label{infbetas}
\cE\{|u_0 -\hat{u}_{0|{\rm past \;\&\; future}}|^2\}=\|1-\sum_{k\neq 0}\beta_ke^{jk\theta}\|_{d\mu}^2.
\end{equation}

In general, the power spectrum $d\mu$ is a bounded nonnegative measure on $[-\pi,\,\pi)$ and admits a decomposition $d\mu=d\mu_\s+fd\theta$ with $d\mu_\s$ a singular measure and $fd\theta$ the absolutely continuous part of $d\mu$ (with respect to the Lebesgue measure). In general the singular part has no effect on the minimal variance of the error, and the corresponding component of $u_k$ can be estimated with arbitrary accuracy using any ``one-sided'' infinite past. The variance of the optimal
one-step-ahead prediction error depends only on the absolutely continuous part of the power spectrum and is given in terms by the celebrated Szeg\"{o}-Kolmogorov formula stated below  (see \cite{Szego} and also, \cite[page 183]{GrenanderSzego}, \cite[Chapter 6]{Varadhan}, \cite{Haykin,StoicaMoses}).\\

\begin{thm}{\sf With $d\mu=d\mu_\s+fd\theta$ as above
\[
\inf_\alpha \|1-\sum_{k>0}\alpha_ke^{jk\theta}\|_{d\mu}^2=\exp\left\{\frac{1}{2\pi}\int_{-\pi}^\pi \log f(\theta)d\theta\right\}
\]
when $\log f\in L_1[-\pi,\pi)$, and zero otherwise.
}\end{thm}
\vspace*{.05in}

In case $\log f\in L_1[-\pi,\pi)$ the prediction-error variance is nonzero and the random process is non-deterministic in the sense of Kolmogorov. In this case, it can be shown that
\[
f(\theta)=\frac{g_f}{|a_f(e^{j\theta})|^2}
\]
where $a_f(z)$ is an outer function in the Hardy space $H_2(\mD)$ with $a_f(0)=1$, i.e.,
\[
a_f(z)=1+a_1z+a_2z^2+\ldots
\]
is analytic in the unit disc $\mD:=\{z\,:\,|z|<1\}$
and its radial limits are square integrable (see \cite{Rudin}).
%\begin{eqnarray*}
%\|a\|_2&:=&\sqrt{\lim_{r\nearrow 1} \intpi |a(re^{j\theta})|^2d\theta}\\
%&=&\sqrt{1+\sum_{k=1}^\infty |a_k|^2} <\infty.
%\end{eqnarray*}
Then, the linear combination
\begin{equation}\label{upast}
\hat{u}_{0|{\rm past}}:=\sum_{k>0} (-a_k) u_{-k}
\end{equation}
serves as the optimal predictor of $u_0$ based on past observations and the least variance of the optimal prediction error becomes
\begin{eqnarray*}
\cE\{|u_0 -\sum_{k>0} (-a_k) u_{-k}|^2\}&& \\
&& \hspace*{-2cm}=\lim_{r\nearrow 1} \intpi |a(re^{j\theta})|^2 f(\theta)d\theta\\ && \hspace*{-2cm}=\exp(\intpi\log(f(\theta))d\theta\\ && \hspace*{-2cm} =: g_f.
\end{eqnarray*}

Analogous expressions exist for the optimal smoothing error and the corresponding smoothing filter which uses both past and future values of $u_\ell$. It is quite interesting, and rather straightforward, that while the variance of the optimal one-step-ahead prediction error is the {\em geometric mean} of the spectral density function, the variance of the error, when a smoothing filter utilizes both past and future, turns out to be the {\em harmonic mean} of the spectral density function.\\

\begin{thm} (see \cite{ansatz}) \label{prop:minvalue}{\sf
With $d\mu=d\mu_\s+fd\theta$ as above
\begin{eqnarray}\nonumber
\inf_{\beta}\|1-\sum_{k\neq 0}\beta_ke^{jk\theta}\|_{d\mu}^2&=&\left(\frac{1}{2\pi}\int_{-\pi}^{\pi} f(\theta)^{-1}d\theta\right)^{-1}\\
&=:& h_f\label{minimalvalue}
\end{eqnarray}
when $f^{-1}\in L_1[-\pi,\pi)$, and zero otherwise.
}
\end{thm}
\vspace*{.05in}

In case $f^{-1}\in L_1[-\pi,\pi)$ the variance of the optimal smoothing error is nonzero and the random process is nondeterministic in the sense that past and future specify the present which can be estimated with zero variance. In this case (see  \cite{ansatz})
\begin{eqnarray*}
b_f(\theta)&=&\ldots +b_2e^{-2j\theta}+b_{-1}e^{-j\theta}+1+\\
&&+ \, b_1e^{j\theta}+b_2e^{2j\theta}+\ldots\\
&=&h_f f(\theta)^{-1}
\end{eqnarray*}
is the image of the optimal smoothing error
$u_0-\sum_{k\neq 0}(-b_k)u_k$
under the Kolmogorov map, and that
\begin{eqnarray*}
\cE\{|u_0-\sum_{k\neq 0}(-b_k)u_k|^2\}&=&\intpi |b_f(\theta)|^2 f(\theta)d\theta\\
&& \hspace*{-1cm}=\intpi h_f^2\left(f(\theta)\right)^{-2}f(\theta)d\theta\\
&& \hspace*{-1cm}=h_f.
\end{eqnarray*}

\section{Degradation of the prediction error variance}
We now consider two distinct spectral density functions $f_1,f_2$ and postulate a situation where filtering of an underlying random process is attempted based on the incorrect choice between these two alternatives.  The variance is then compared with the least possible variance which is achieved when the correct choice is made (i.e., when the predictor is optimal for the spectral density against which it is being evaluated). The degradation of performance is quantified by how much the ratio of the two prediction-error variances exceeds the identity.
This ratio serves as a measure of mismatch between the two spectral densities (the one which was used to design the predictor and the one against which it is being evaluated).
The resulting mismatch turns out to be {\em scale-invariant} ---i.e., the expression is homogeneous. Hence, as a measure of distance it actually quantifies distance between the positive rays that the two spectral density functions define, and thus, it quantifies distance between the respective ``shapes.'' It turns out that this distance is convex on logarithmic intervals and has a number of distance-like properties, short of being a metric.

Let us assume that both $\log f_1, \log f_2 \in L_1[-\pi,\pi)$ and hence, that
\[
f_i(\theta)=\frac{g_{f_i}}{|a_{f_i}(e^{j\theta})|^2}
\]
for corresponding outer $H_2$-functions $a_{f_i}(z)$ normalized as before so that $a_{f_i}(0)=1$, for $i\in\{1,2\}$. Obviously,
\[
g_{f_i}=\exp\left(\intpi \log(f_i(\theta))d\theta\right)
\]
denotes the geometric mean of $f_i$ as before, for $i\in\{1,2\}$. These expressions represent the least variances when the predictor is chosen on the basis of the correct spectral density function. If however, the predictor is based on $f_2$ whereas the underlying process has $f_1$ as its spectral density, then the variance of the prediction error turns out to be
\begin{eqnarray*}
\intpi |a_{f_2}(e^{j\theta})|^2 f_1(\theta)d\theta &=& \left(\intpi \frac{f_1(\theta)}{f_2(\theta)}d\theta\right)g_{f_2}.
\end{eqnarray*}
If we divide this variance by the optimal value $g_{f_1}$ we obtain
\begin{eqnarray}
{\rho}_{a/g}(f_1,f_2)&:=&
\left(\intpi \frac{f_1(\theta)}{f_2(\theta)}d\theta\right)\frac{g_{f_2}}{g_{f_1}}\\\nonumber
&&\hspace*{-3cm}=
\left(\intpi \frac{f_1(\theta)}{f_2(\theta)}d\theta\right)\frac{\exp(\intpi \log f_2(\theta) d\theta)}{\exp(\intpi \log f_1(\theta) d\theta)}\\
&&\hspace*{-3cm}=
\left(\intpi \frac{f_1(\theta)}{f_2(\theta)}d\theta\right)\frac{1}{\exp\left(\intpi \log(\frac{f_1(\theta)}{f_2(\theta)})d\theta\right)}.
\end{eqnarray}
This is the ratio of the {\em arithmetic mean} over the {\em geometric mean} of the fraction $f_1/f_2$ of the two spectral density functions. The expression ${\rho}_{a/g}(f_1,f_2)$ is not symmetric in the two arguments. The subscript ``$_{a/g}$'' signifies ratio of arithmetic over geometric means.

The logarithm of ${\rho}_{a/g}(f_1,f_2)$ is nonnegative and defines a notion of distance between rays of density functions. Henceforth, we denote this logarithm by
\begin{eqnarray}
\label{deltaag}
{\delta}_{a/g}(f_1,f_2)&:=&\log\left( {\rho}_{a/g}(f_1,f_2)\right)\\
\nonumber &=& \log\left(\intpi \frac{f_1(\theta)}{f_2(\theta)}d\theta\right)
%\\&&
-\intpi \log\left(\frac{f_1(\theta)}{f_2(\theta)}\right)d\theta.\nonumber
\end{eqnarray}
Alternatively, we can view the above as slackness of a Jensen-type inequality.

Before we discuss key properties of ${\delta}_{a/g}$, we introduce a natural class of paths connecting density functions: for any two density functions $f_a,f_b$,
\[
f_{\tau,\,a,b}:=f_a^{1-\tau}f_b^\tau, \mbox{ for }\tau\in[0,1],
\]
defines a {\em logarithmic interval} between $f_a$ and $f_b$. The terminology stems from the fact that whenever the needed logarithms exist,
\[
f_{\tau,\,a,b}=e^{(1-\tau)\log(f_a)+\tau \log(f_b)}, \mbox{ for }\tau\in[0,1].
\]
Later on we will see that these represent geodesics on the manifold of density functions with respect to an induced pseudo-Riemannian metric.

\begin{prop}\label{prop:prop3}{\sf Let $f_i$, $i\in\{1,2,3\}$ represent density functions defined on $[-\pi,\pi)$. The following hold:\\
\begin{tabular}{rl}
(i)&${\delta}_{a/g}(f_1,f_2)\in \mR_+\cup \{\infty\}$.\\
(ii)&${\delta}_{a/g}(f_1,f_2)=0$ $\Leftrightarrow$ $f_1(\theta)/f_2(\theta)$ is constant.\\
(iii)&$\delta_{a/g}(f_1,f_{\tau,\,1,b})$ is monotonically increasing\\
& for $\tau\in[0,1]$ and $b\in\{2,3\}$.\\
(iv)&$\delta_{a/g}(f_1,f_{\tau,\,2,3})$ is convex in $\tau$.
%\item[(iii)] ${\delta}_{ag}(f_1,f_2)+{\delta}_{ag}(f_2,f_3)\geq {\delta}_{ag}(f_1,f_3)$.
\end{tabular}
}\end{prop}
\vspace*{.05in}

\begin{proof} Claims (i-ii) follow from the fact that the arithmetic mean of a function always exceeds the geometric mean, and that they are equal whenever a function is constant. In particular, the ordering, as to which is larger, follows from Jensen's inequality
\[ \log\left(\intpi f(\theta)d\theta\right) \geq \intpi \log\left( f(\theta)\right)d\theta
\]
for any $f(\theta)\geq 0$. The fact that the two are equal only when $f(\theta)$ is constant can be obtained easily using a variational argument. Then (i-ii) follow, when we substitute $f=f_1/f_2$ and then take the logarithm.

Next we show (iv) and use it to deduce (iii). Since
\begin{eqnarray*}
\delta(f_1,f_a^{(1-\tau)}f_b^\tau)&=&\log\intpi \frac{f_1(\theta)}{f_a^{(1-\tau)}(\theta)f_b^\tau(\theta)}d\theta\\
&&-
\intpi \log\left(\frac{f_1(\theta)}{f_a^{(1-\tau)}(\theta)f_b^\tau(\theta)}\right) d\theta\\
&=& \log \intpi \left(\frac{f_1(\theta)}{f_a(\theta)}\right)\left(\frac{f_a(\theta)}{f_b(\theta)}\right)^\tau d\theta\\
&&-
\intpi \log \left(\frac{f_1(\theta)}{f_a(\theta)}\right) d\theta\\&&-
\intpi \tau \log\left(\frac{f_a(\theta)}{f_b(\theta)}\right) d\theta,
\end{eqnarray*}
the derivative with respect to $\tau$ becomes
\begin{eqnarray*}
\frac{d}{d\tau}\delta(f_1,f_a^{(1-\tau)}f_b^\tau)&=&\frac{1}{\int_{-\pi}^\pi  \frac{f_1(\theta)}{f_a(\theta)}\left(\frac{f_a(\theta)}{f_b(\theta)}\right)^\tau d\theta}\\&&\hspace*{-2cm}\times
\int_{-\pi}^\pi  \frac{f_1(\theta)}{f_a(\theta)}\left(\log\frac{f_a(\theta)}{f_b(\theta)}\right)\left(\frac{f_a(\theta)}{f_b(\theta)}\right)^\tau d\theta
\\&&\hspace*{-2cm}-\intpi \log\left(\frac{f_a(\theta)}{f_b(\theta)}\right) d\theta
\end{eqnarray*}
and the second derivative,
\begin{eqnarray*}
\frac{d^2}{d\tau^2}\delta(f_1,f_a^{(1-\tau)}f_b^\tau)&=&
\frac{-1}{\left(\int_{-\pi}^\pi  \frac{f_1(\theta)}{f_a(\theta)}\left(\frac{f_a(\theta)}{f_b(\theta)}\right)^\tau d\theta\right)^2}
\\&&\hspace*{-2cm}\times\left(\int_{-\pi}^\pi  \frac{f_1(\theta)}{f_a(\theta)}\left(\log\frac{f_a(\theta)}{f_b(\theta)}\right)\left(\frac{f_a(\theta)}{f_b(\theta)}\right)^\tau d\theta \right)^2\\
&&\hspace*{-4cm}+
\frac{1}{\int_{-\pi}^\pi  \frac{f_1(\theta)}{f_a(\theta)}\left(\frac{f_a(\theta)}{f_b(\theta)}\right)^\tau d\theta}
\int_{-\pi}^\pi  \frac{f_1(\theta)}{f_a(\theta)}\left(\log\frac{f_a(\theta)}{f_b(\theta)}\right)^2\left(\frac{f_a(\theta)}{f_b(\theta)}\right)^\tau d\theta
\end{eqnarray*}
But from Cauchy's inequality we have that
\begin{eqnarray*}
\left(\int_{-\pi}^\pi  \frac{f_1(\theta)}{f_a(\theta)}\left(\frac{f_a(\theta)}{f_b(\theta)}\right)^\tau d\theta\right) &&\\
&&\hspace*{-2.7cm}\times
\left(\int_{-\pi}^\pi  \frac{f_1(\theta)}{f_a(\theta)}\left(\log\frac{f_a(\theta)}{f_b(\theta)}\right)^2 \left(\frac{f_a(\theta)}{f_b(\theta)}\right)^\tau d\theta\right)\\
&&\hspace*{-2.7cm}
\geq \left(\int_{-\pi}^\pi  \frac{f_1(\theta)}{f_a(\theta)}\left(\log\frac{f_a(\theta)}{f_b(\theta)}\right) \left(\frac{f_a(\theta)}{f_b(\theta)}\right)^\tau d\theta\right)^2.
\end{eqnarray*}
Hence, the second derivative is nonnegative. Claim (iv) is seen to hold true after we set $a=2$ and $b=3$. To establish claim (iii) set $a=1$ and $b=2$ in the above. Then,
\begin{eqnarray*}
\frac{d^2}{d\tau^2}\delta(f_1,f_1^{(1-\tau)}f_2^\tau)\geq 0.
\end{eqnarray*}
But $\delta(f_1,f_1^{(1-\tau)}f_2^\tau)=0$ for $\tau=0$, and $\delta(f_1,f_1^{(1-\tau)}f_2^\tau)\geq 0$ for all $\tau$. Hence,
the derivative at $\tau=0$ must be nonnegative and $\delta(f_1,f_1^{(1-\tau)}f_1^\tau)$ must increase as $\tau\nearrow 1$. This completes the proof of (iii), and the proof of the proposition.
\end{proof}

Since $\delta_{a/g}(f_1,f_2)$ is not symmetric in its arguments, it is quite natural to consider the symmetrized version
\[
\delta(f_1,f_2):=\delta_{a/g}(f_1,f_2)+\delta_{a/g}(f_2,f_1).
\]
All properties listed in Proposition \ref{prop:prop3} hold true for $\delta(f_1,f_2)$ as well. Furthermore, interestingly,
\begin{eqnarray*}
\delta(f_1,f_2)&=& \log\left(\intpi \frac{f_1(\theta)}{f_2(\theta)}d\theta\right) + \log\left(\intpi \frac{f_2(\theta)}{f_1(\theta)}d\theta\right)\\
&=&\log\left(\left(\intpi \frac{f_1(\theta)}{f_2(\theta)}d\theta\right)\left(\intpi \frac{f_2(\theta)}{f_1(\theta)}d\theta\right)\right)
\end{eqnarray*}
which is the {\em logarithm of the ratio of the arithmetic mean over the harmonic mean} of the ``likelihood'' fraction $f_1/f_2$. Again, the distance of this ratio from one quantifies how far $f_1/f_2$ is from being constant. We now summarizing  the claimed properties of $\delta(\cdot,\cdot)$.

\begin{prop}\label{prop:prop4}{\sf Let $f_i$, $i\in\{1,2,3\}$, density functions on $[-\pi,\pi)$. The following hold:\\
\begin{tabular}{rl}
(i)&${\delta}(f_1,f_2)\in \mR_+\cup \{\infty\}$.\\
(ii)&${\delta}(f_1,f_2)=0$ $\Leftrightarrow$ $f_1(\theta)/f_2(\theta)$ is constant.\\
(iii)&$\delta(f_1,f_{\tau,\,1,b})$ is monotonically increasing\\ &for $\tau\in[0,1]$ and $b\in\{2,3\}$.\\
(iv)&$\delta(f_1,f_{\tau,\,2,3})$ is convex in $\tau$.\\
%\item[(iii)] ${\delta}_{ag}(f_1,f_2)+{\delta}_{ag}(f_2,f_3)\geq {\delta}_{ag}(f_1,f_3)$.
\end{tabular}
}\end{prop}
\vspace*{.05in}

\begin{proof} Properties (i), (ii), and (iv) are a direct consequence of the corresponding properties given in Proposition \ref{prop:prop3} for $\delta_{a/g}(\cdot,\cdot)$. Property (iii) on the other hand follows as before from (iv) and the fact that the derivative of $\delta(f_1,f_{\tau,\,1,b})$ at $\tau=0$ is zero.
\end{proof}

\section{An example}

\begin{figure}[htb]
\begin{center}
\psfrag{f_x}{$c_1$}
\includegraphics[totalheight=7cm]
{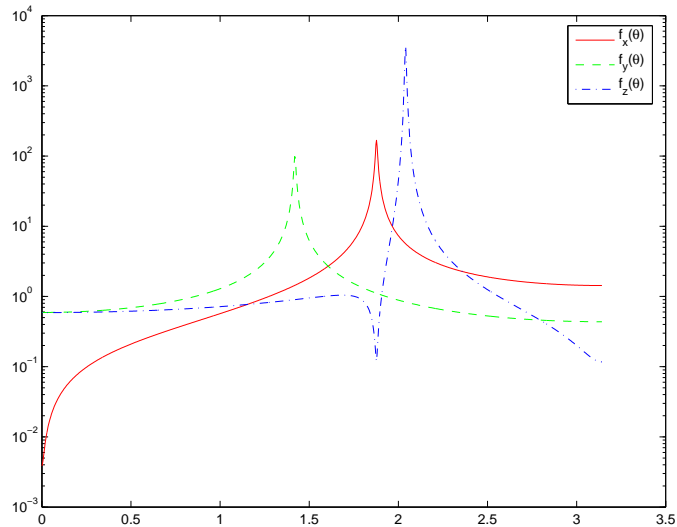}
\caption{Three power spectral densities $f_1$ ({\color{red}--}), $f_2$ ({\color{green}- -}), $f_3$ ({\color{blue}-$\cdot$-})}
\label{fig:spectra}
\end{center} 
\end{figure}

\begin{figure}[htb]
\begin{center}
\includegraphics[totalheight=7cm]
{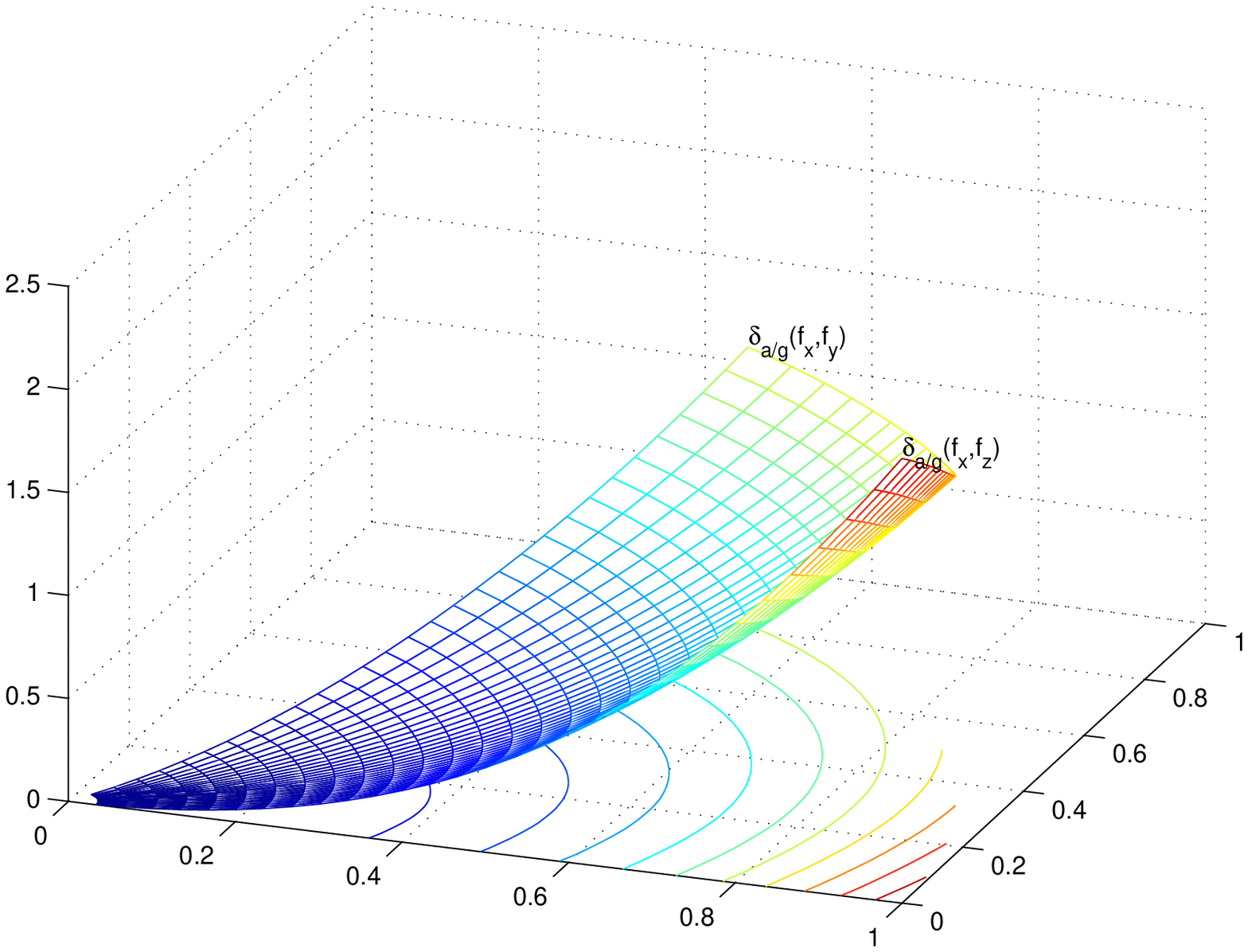}
\caption{$\delta_{a/g}(f_1,f_1^{1-\tau} f_2^{\tau-\sigma\tau}f_3^{\sigma\tau})$}
\label{fig:surf1}
\end{center} 
\end{figure}

\begin{figure}[htb]
\begin{center}
\includegraphics[totalheight=7cm]
{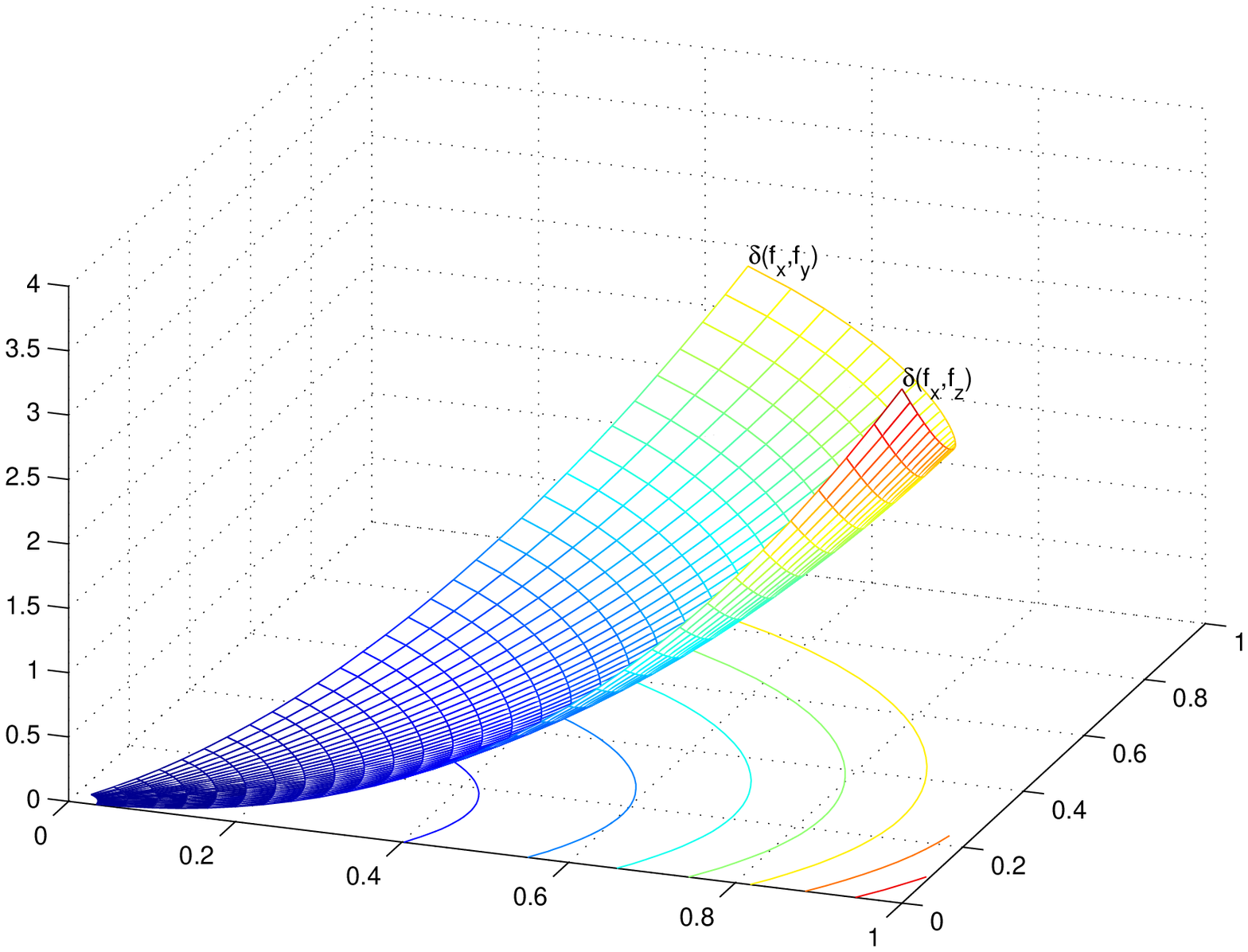}
\caption{$\delta(f_1,f_1^{1-\tau} f_2^{\tau-\sigma\tau}f_3^{\sigma\tau})$}
\label{fig:surf2}
\end{center} 
\end{figure}

\begin{figure}[htb]
\begin{center}
\includegraphics[totalheight=7cm]
{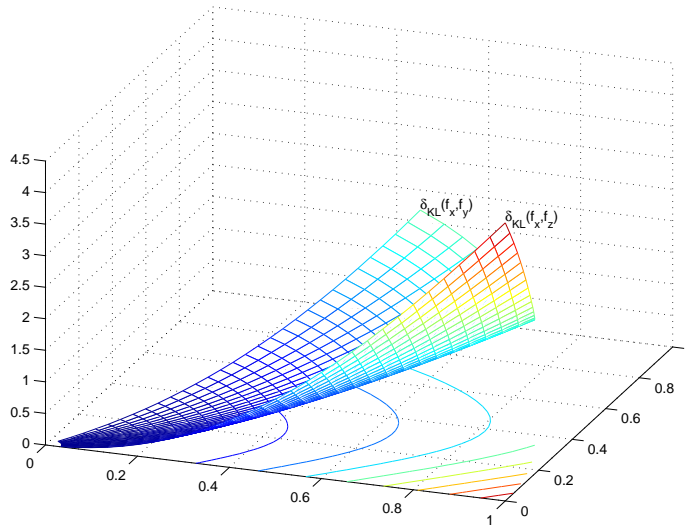}
\caption{$\delta_{\rm KL}(f_1,f_1^{1-\tau} f_2^{\tau-\sigma\tau}f_3^{\sigma\tau})$}
\label{fig:surf3}
\end{center} 
\end{figure}

In order to illustrate the quantitative behavior of these measures,
we consider three specific power spectra labeled $f_1,f_2,f_3$, as before. These are shown in Figure \ref{fig:spectra}. We then consider the triangle formed with those power spectra as vertices and connected using logarithmic intervals. The interior of the triangle is similarly sampled at logarithmically placed points. In essence, we consider the family of power spectral densities
\begin{equation}\label{points}
f_1^{(1-\tau)} \left( f_2^{(1-\sigma)}f_3^\sigma\right)^\tau \mbox{ for }\tau,\sigma\in[0,1].
\end{equation}
For each value of $\tau,\sigma$ (sampled appropriately), we evaluate
$\delta_{a/g}(f_1,f)$, $\delta(f_1,f)$ and compare these to the Kullback-Leibler divergence between the suitably normalized functions
\[
\hat{f}_k(\theta)=\frac{f_k(\theta)}{{\intpi f_k d\theta}},\;\; k\in\{1,2,3\}.
\]
The normalization is necessary if the Kullback-Leibler divergence is to have properties of a distance measure (i.e., nonnegative when its arguments are different, etc.).
Thus, we denote
\begin{eqnarray}
\delta_{\rm KL}(f_1,f_2)\nonumber
&:=&\intpi \hat{f}_1(\theta)\log\left(\frac{\hat{f}_1(\theta)}{\hat{f}_2(\theta)}\right)d\theta\\
&=&
\frac{1}{\intpi f_1(\theta)d\theta}\left( \intpi f_1(\theta) \left(\log(\frac{f_1(\theta)}{f_2(\theta)})
%\right.\\&&\left.
- \log\frac{\intpi f_1(\theta)d\theta}{\intpi f_2(\theta) d\theta}\right)d\theta\right).\label{eq:KL}
\end{eqnarray}
The set of power spectra in (\ref{points}) is thought of as a set of points forming an equilateral triangle, conceptually sitting on the $xy$-plane. Then, the vertical axis represents distance from $f_1$, measured using these three alternative measures. The corresponding surfaces are drawn in Figures \ref{fig:surf1}-\ref{fig:surf3}.
%We also draw the $L_2$-distance from $f_1$ for comparison in Figure  \ref{fig:surfL2}.
The three power spectra used are as follows:
\begin{eqnarray*}
f_1(\theta) &=& \left|\frac{(z-.99)}{(z^2+.6z+.99)}\right|^2_{z=e^{j\theta}}\\
f_2(\theta)  &=& \left|\frac{1}{(z^2-.3z+.99)}\right|^2_{z=e^{j\theta}}\\
f_3(\theta)  &=& \left|\frac{(z+.9)(z^2+.6z+.99)}{(z^2+.9z+.99)(z^2+.9z+.99)}\right|^2_{z=e^{j\theta}}
\end{eqnarray*}

%\begin{figure}[htb]
%\begin{center}
%\includegraphics[totalheight=5.2cm]
%{/Users/tryphon/almostmetrics/surfaceL2.eps}
%\caption{$\sqrt{\intpi (f_1-f_1^{1-\tau} f_2^{\tau-\sigma\tau}f_3^{\sigma\tau})^2d\theta}$}
%\label{fig:surfL2}
%\end{center} 
%\end{figure}

There appears to be little qualitative difference between $\delta(\cdot,\cdot)$,  $\delta_{a/g}(\cdot,\cdot)$, and $\delta_{\rm KL}(\cdot,\cdot)$.
They are also quite similar in that it is easy to calculate functional forms for minimizers of either of these distance measures under moment constraints (see Section \ref{sec:minimizers} below).
Hence, it is important to undercore that $\delta_{\rm KL}(\cdot,\cdot)$ lacks an intrinsic interpretation as a distance measure between power spectra, in contrast to $\delta(\cdot,\cdot)$ which therefore may be preferable for exactly that reason.

\section{Functional form of minimizers in moment problems}\label{sec:minimizers}

A large class of spectral analysis problems is typified by the trigonometric moment problem where a power spectral density is sought to match a partial sequence of autocorrelation samples, i.e., a positive function $f$ is sought such that
\begin{equation}\label{momentconstraints}
R_k=\intpi e^{jk\theta}f(\theta)d\theta,\mbox{ for }k\in\{0,1,\ldots,n\},
\end{equation}
see e.g., \cite{GrenanderSzego,csiszar2,Haykin,BGL2,GL}.
Since, in general, the family of consistent $f$'s is large, a particular one is chosen ``closest'' to a given ``prior''. Maximum entropy spectral analysis, for instance, can be interpreted as seeking the spectral density closest in the Kullback-Leibler sense to one which is flat, i.e., the prior in this case is the power spectral density of white noise (see e.g., \cite{GL}). In the same spirit we may pose the problem of seeking $f$ closest to $f_{\rm prior}$ in the sense of minimizing e.g., $\delta_{a/g}(f,f_{\rm prior})$ and subject to the moment constraints (\ref{momentconstraints}).

To this end, as usual, we introduce Lagrange multipliers $\lambda_k$ ($k\in\{0,1,\ldots,n\}$) and form the Lagrangian
\[
\cL(f,\lambda_0,\ldots,\lambda_n):=\log\intpi \frac{f(\theta)}{f_{\rm prior}(\theta)}d\theta
-\intpi \log\frac{f(\theta)}{f_{\rm prior}(\theta)}d\theta +\sum_{k=-n}^n
\lambda_k\left(R_k-\intpi e^{jk\theta}f(\theta)d\theta\right).
\]
Setting the variation of $\cL$ identically to zero for all perturbations of $f$ gives conditions that help identify the functional form of minimizing $f$'s. Briefly,
\begin{eqnarray*}
\delta\cL&=& \log \intpi \frac{f+\Delta}{f_{\rm prior}}d\theta -\intpi \log\frac{f+\Delta}{f_{\rm prior}}d\theta\\
&&-\log \intpi \frac{f}{f_{\rm prior}}d\theta+\intpi \log\frac{f}{f_{\rm prior}}d\theta - \sum_{k=-n}^n
\lambda_k\left(\intpi e^{jk\theta}\Delta(\theta)d\theta\right)\\
&\simeq& \frac{ \intpi \frac{\Delta}{f_{\rm prior}}d\theta}{ \intpi \frac{f}{f_{\rm prior}}d\theta}
- \intpi \frac{\Delta}{f}d\theta -\sum_{k=-n}^n
\lambda_k\left(\intpi e^{jk\theta}\Delta d\theta\right)\end{eqnarray*}
after we eliminate higher order terms. Stationarity conditions require
that the above is identically zero for all (small) functions $\Delta$. This leads to
\[
\frac{1}{f_{\rm prior}(\theta){ \intpi {f(\theta)/f_{\rm prior}(\theta)}d\theta}}
-  \frac{1}{f(\theta)} -\sum_{k=-n}^n
\lambda_k e^{jk\theta}  = 0
\]
from which we deduce that a minimizing $f$ must be of the form
\begin{eqnarray}\label{eq:functionalform}
f(\theta)= \frac{\kappa f_{\rm prior}(\theta)}{1 - \kappa f_{\rm prior}(\theta)\sum_{k=-n}^n
\lambda_k e^{jk\theta}}
\end{eqnarray}
with
\begin{equation}\label{extra}
\kappa=\intpi f(\theta)/f_{\rm prior}(\theta)d\theta.
\end{equation}
Then, values for $\kappa$ as well as for the Lagrange multipliers must be determined so that $f$ in (\ref{eq:functionalform}) satisfies (\ref{momentconstraints}) and (\ref{extra}) ---this can be done for instance using homotopy methods in \cite{ac_may2004,IT}. It is interesting that when $f_{\rm prior}=1$, the minimizer is the same as in the one obtained by applying the maximum entropy principle (e.g., see \cite{GL,Haykin}), i.e., it turns out to be an all-pole spectral density function which is of course uniquely identified by the moment constraints. Evidently, in general, minimizing $\delta_{a/g}(f,f_{\rm prior})$ gives a different answer than the one obtained by minimizing $\delta_{\rm KL}(f,f_{\rm prior})$ or, by minimizing other distances. Yet, all such problems are similar and can be dealt with in two steps. First identify the functional form of a minimizer and then determine values for the coefficients so as to satisfy (\ref{momentconstraints}). The latter step requires solving a nonlinear problem in general, and can be approached in a variety of ways (e.g., as in \cite{ac_may2004,IT}).

\section{Riemannian metrics and geodesics}\label{sec:riemannian}

Infinitesimal perturbations about a given power spectral density function, when measured by any of $\delta_{a/g}(\cdot,\cdot)$,  $\delta(\cdot,\cdot)$,  or $\delta_{\rm KL}(\cdot,\cdot)$, give rise to nonnegative definite quadratic forms.
These forms are in fact nonsingular on directions other than rays emanating from the origin. This is due to the fact that the aforementioned distances do not separate points on such rays while they give nonzero distance otherwise. They thus induce Riemannian metrics on suitably defined manifold of spectral rays. In this section (and the current paper) we focus on the particular metric induced by $\delta(\cdot,\cdot)$, we show how to characterize geodesics, and verify that logarithmic intervals are in fact geodesics.

Throughout we assume that all functions are smooth enough so that the indicated integrals exist. This, in particular, can be ensured if all spectral density functions are bounded and have bounded inverses as well as bounded derivatives. Weaker conditions are clearly possible. For the purposes of this section we define
\begin{eqnarray*}
\cF&:=&\{f \;:\; f(\theta) \mbox{ differentiable on }[-\pi,\pi],\\
&& \mbox{ with  }f(\theta)> 0, \mbox{ and both }f(\theta),\frac{df(\theta)}{d\theta} \mbox{ square integrable}\}.
\end{eqnarray*}
With a suitable norm on $f$ and its derivative, $\cF$ becomes a (Banach) manifold.
We also recall the definition
$%\begin{eqnarray*}
\|f\|_k:=\left(\intpi |f(\theta)|d\theta\right)^{1/k}
$%\end{eqnarray*}
of the $k$-th norm, applicable to any $f$ on $[0,2\pi]$
provided the integral exist. Whenever, $f$ is a density function, the absolute value is obsiously unnecessary.

\begin{prop}\label{prop:deltas}{\sf
Let $f,f+\Delta\in\cF$ where $\Delta$ is a perturbation such that $|\Delta/f|<1$. Then
\begin{eqnarray*}
\delta(f,f+\Delta) &=&\intpi\left(\frac{\Delta(\theta)}{f(\theta)}\right)^2d\theta- 
\left(\intpi\frac{\Delta(\theta)}{f(\theta)}d\theta\right)^2 + O(\|\Delta/f\|^3_1)\\
\delta_{\rm a/g}(f,f+\Delta) &=&\frac12\left(\intpi\left(\frac{\Delta(\theta)}{f(\theta)}\right)^2d\theta- 
\left(\intpi\frac{\Delta(\theta)}{f(\theta)}d\theta\right)^2\right) + O(\|\Delta/f\|^3_1)\\
\delta_{\rm KL}(f,f+\Delta) &=&\frac12\left(\intpi\frac{\left(\Delta(\theta)\right)^2}{f(\theta)}d\theta- 
\left(\intpi\Delta(\theta)d\theta\right)^2\right) + O(\|\Delta/f\|^3_1).
\end{eqnarray*}
Here, $O(\|\Delta/f\|^3_1)$ indicates terms of order $3$ or higher.
}\end{prop}
\vspace*{.05in}

\begin{proof}
We prove only the first claim. The other two can be shown in an identical manner.
We expand
$\delta(f,f+\Delta)$
using the series $\log(1+x)=x-\frac{1}{2}x^2+\frac{1}{3}x^3-\ldots$, as follows
\begin{eqnarray*}
\delta(f,f+\Delta)&=&
 \log\left(\intpi\frac{f(\theta)+\Delta(\theta)}{f(\theta)}d\theta\right)+
 \log\left(\intpi\frac{f(\theta)}{f(\theta)+\Delta(\theta)}d\theta\right)\\
 &=&\log\left(1+\intpi\frac{\Delta(\theta)}{f(\theta)}d\theta\right)+
 \log\left(\intpi\frac{1}{1+\frac{\Delta(\theta)}{f(\theta)}}d\theta\right)\\
 &=& \intpi\frac{\Delta(\theta)}{f(\theta)}d\theta - \frac{1}{2} \left(\intpi\frac{\Delta(\theta)}{f(\theta)}d\theta\right)^2 + O(\|\frac{\Delta(\theta)}{f(\theta)}\|^3_1)\\
 &&+ \log\left(
 \intpi(1-\frac{\Delta(\theta)}{f(\theta)}+\left(\frac{\Delta(\theta)}{f(\theta)}\right)^2 + \ldots)d\theta\right)\\
 &=& \intpi\frac{\Delta(\theta)}{f(\theta)}d\theta - \frac{1}{2} \left(\intpi\frac{\Delta(\theta)}{f(\theta)}d\theta\right)^2  \\
&&+\intpi (-\frac{\Delta(\theta)}{f(\theta)} + \left(\frac{\Delta(\theta)}{f(\theta)}\right)^2) d\theta
-\frac{1}{2}\left(\intpi(-\frac{\Delta(\theta)}{f(\theta)})d\theta\right)^2
 + O(\|\frac{\Delta}{f}\|^3_1)
\end{eqnarray*}
which proves the first claim after canceling and collecting terms. The other two expressions can be shown similarly. Note that $\|\Delta/f\|_k\leq \|\Delta/f\|_1$ for all $k$, since $|\Delta/f|<1$.
\end{proof}

Since $\delta_{a/g}(\cdot,\cdot)$,  $\delta(\cdot,\cdot)$,  and $\delta_{\rm KL}(\cdot,\cdot)$ do not separate power spectra which are scalar multiple of one another, we may consider equivalence classes
\begin{eqnarray*}
(f)_\sim &:=&\{f_1\in\cF \;:\; f_1=cf,\,c\in\mR_+\}.
\end{eqnarray*}
obtained from any spectral density function $f$ by a scaling factor, constant across $[-\pi,\pi]$. These can be thought of as rays. They can be identified by pointing to one particular representative. Thus, in particular, the set of rays can be identified with
\[
\fF:=\{f\;:\;f\in\cF,\; \|f\|_1=1\}.
\]
This set can be given the structure of a manifold and thought of as a set of probability density functions on $[-\pi,\pi]$. Alternatively, we can consider spectral densities as belonging to $\cF$ and accept the fact that we have a pseudo-Riemannian metric which vanishes along certain directions.
This follows from Proposition \ref{prop:deltas}, since $\delta(\cdot,\cdot)$ (and similarly, $\delta_{a/g}$ and $\delta_{\rm KL}$) defines a nonnegative definite quadratic form\footnote{
If we mod-out scaling and stay on a manifold of ``spectral rays'', this quadratic form becomes positive definite, and defines a Riemannian metric.} at each ``point'' $f\in\cF$ via
\begin{equation}\label{Riemannianmetric}
\Delta \mapsto \intpi\left(\frac{\Delta(\theta)}{f(\theta)}\right)^2d\theta- 
\left(\intpi\frac{\Delta(\theta)}{f(\theta)}d\theta\right)^2.
\end{equation}

A rather standard way to measure distances on manifold is to trace geodesics connecting points, and compute the length of such paths. Thus, it may be of interest to characterize geodesics in our case as well. We refrain from excessively technical jargon and, following the earlier suggestion, simply assume that all integrals exist.

Consider a path $f_\tau$, $\tau\in[0,1]$, of spectral density functions connecting two given ones, namely $f_0$ and $f_1$. Note that $f_\tau$ is a function of two arguments,
the path parameter $\tau$ and the frequency $\theta$ ---hence, we often write $f_\tau(\theta)$. The length traversed as $\tau$ varies from $0$ to $1$ is simply
\begin{eqnarray}\label{eq:length}
\ell(f_\tau:0,1)&:=&\int_0^1 \sqrt{\delta(f_\tau,f_{\tau+d\tau})}\\
&=&\int_0^1\sqrt{\log\left( \intpi \frac{f_\tau(\theta)+\dot{f}_\tau(\theta)d\tau}{f_\tau(\theta)} d\theta \right)
+\log\left( \intpi \frac{f_\tau(\theta)}{f_\tau(\theta)+\dot{f}_\tau(\theta)d\tau} d\theta \right)}\nonumber\\
&=&
\int_0^1\sqrt{
\intpi\left(\frac{\dot{f}_\tau(\theta)}{f_\tau(\theta)}\right)^2d\theta- 
\left(\intpi\frac{\dot{f}_\tau(\theta)}{f_\tau(\theta)}d\theta\right)^2}d\tau. \label{eq:interesting}
\end{eqnarray}
In the last step we eliminated higher order terms in $d\tau$ inside the integral, since those integrate to zero.
Here and throughout ``$\;\dot{\;}\;$'' (dot), as in $\dot{f}$ is used to denote derivative with respect to $\tau$, i.e.,
\[
\dot{f_\tau}(\theta):=\frac{\partial{f_\tau}(\theta)}{\partial \tau}.
\]
Interestingly, the expression in (\ref{eq:interesting}) only depends on $\dot{f}_\tau(\theta)/f_\tau(\theta)$. Thus, if we define
\[
x_\tau:=\log(f_\tau),
\]
then $\dot{x}_\tau(\theta)=\dot{f}_\tau(\theta)/f_\tau(\theta)$ and
\[
\]
\begin{eqnarray}\label{eq:length2}
\ell(f_\tau:0,1)&:=&\int_0^1\sqrt{
\intpi\left(\dot{x}_\tau(\theta)\right)^2d\theta- 
\left(\intpi\dot{x}_\tau(\theta)d\theta\right)^2}d\tau.
\end{eqnarray}
The requirement that the end point of $f_\tau$ coincide with $f_0$ and $f_1$, readlily translates into boundary conditions for $x_\tau$, namely $x_0=\log f_0$ and $x_1=\log f_1$. The task of finding extremals of such integrals leads to Euler-Lagrange equations for the path $x_\tau$. More specifically, the Lagrangian corresponding to (\ref{eq:length2}) is
\[
L(x_\tau,\dot{x}_\tau,\tau):=\sqrt{
\intpi\left(\dot{x}_\tau(\theta)\right)^2d\theta- 
\left(\intpi\dot{x}_\tau(\theta)d\theta\right)^2}
\]
and only depends on $\dot{x}_\tau$. Therefore $\partial L/\partial x_\tau=0$ and the Euler-Lagrange equations
\[
\frac{\partial L}{\partial x_\tau}-\frac{\partial}{\partial\tau}\frac{\partial L}{\partial \dot{x}_\tau}=0
\]
simplify to $\partial L/\partial \dot{x}_\tau$ being independent of $\tau$. Since $\dot{x}_\tau$ enters in $L$ through an integral over $\theta$, the partial derivative with respect to $\dot{x}_\tau$ is infinitesimal. Thus, we write
\[\frac{\partial L}{\partial \dot{x}_\tau} = v(\theta)d\theta,
\]
which is independent of $\tau$ as we just explained. Since
\begin{eqnarray*}
\frac{\partial L}{\partial \dot{x}_\tau}&=&\frac{1}{2}\frac{1}{\sqrt{
\intpi\left(\dot{x}_\tau(\theta)\right)^2d\theta- 
\left(\intpi\dot{x}_\tau(\theta)d\theta\right)^2}}
\frac{\partial}{\partial\dot{x}_\tau}\left(\intpi\left(\dot{x}_\tau(\theta)\right)^2d\theta- 
\left(\intpi\dot{x}_\tau(\theta)d\theta\right)^2\right)
\end{eqnarray*}
where the latter term produces again a differential in $\theta$, it follows that
\begin{eqnarray*}
2\sqrt{
\intpi\left(\dot{x}_\tau(\theta)\right)^2d\theta- 
\left(\intpi\dot{x}_\tau(\theta)d\theta\right)^2}v(\theta)d\theta&=&
2\dot{x}_\tau(\theta)- 
2\left(\intpi\dot{x}_\tau(\theta)d\theta\right)d\theta.
\end{eqnarray*}
Alternatively,
\begin{eqnarray}\label{eq:must}
\frac{\dot{x}_\tau(\theta)- 
\left(\intpi\dot{x}_\tau(\theta)d\theta\right)}{
\sqrt{
\intpi\left(\dot{x}_\tau(\theta)\right)^2d\theta- 
\left(\intpi\dot{x}_\tau(\theta)d\theta\right)^2}} &=&v(\theta),
\end{eqnarray}
which simply says that the variation of $\dot{x}_\tau$ about the mean, as a function of $\theta$, normalized by a ``standard deviation''-like quantity must be independent of $\tau$. We summarize our conclusion as follows.

\begin{prop}\label{prop:geodesics}{\sf
Given two spectral density functions $f_0,f_1$, extremal (geodesic) paths $f_\tau$ ($\tau\in[0,1]$) connecting the two, in the sense of achieving a local extremal of the path integral $\ell(f_\tau:0,1)$, must satisfy (\ref{eq:must}) for $x_\tau=\log f_\tau$, i.e., the left hand side of (\ref{eq:must}) must be independent of $\tau$.
}\end{prop}
\vspace*{.05in}

\begin{proof} The proof has been established in the arguments leading to the proposition.\end{proof}

We finally verify that logarithmic intervals satisfy (\ref{eq:must}). This is rather straightforward since, for
\[
f_\tau(\theta)=f_0(\theta)\left(\frac{f_1(\theta)}{f_0(\theta)}\right)^\tau,
\]
the logarithm
\begin{eqnarray*}
x_\tau(\theta)&=&\log f_\tau(\theta)\\
&=& \log f_0(\theta) +\tau(\log f_1(\theta)-\log f_0(\theta))
\end{eqnarray*}
is a linear function of $\tau$ and the derivative
\begin{eqnarray*}
\dot{x}_\tau(\theta)&=&\frac{\partial}{\partial\tau} \log(f_\tau(\theta))\\
&=&\log(f_1(\theta))-\log(f_0(\theta))\\
&=&x_1(\theta)-x_2(\theta)
\end{eqnarray*}
is already independent of $\tau$. The ratio $f_1/f_0$ plays a r\^{o}le analogous to the likelihood ratio of probability theory. The length of logarithmic intervals can be computed in terms of this ratio by simple inspection since from (\ref{eq:interesting})
\[
\dot{x}_\tau(\theta)=\frac{\dot{f}_\tau(\theta)}{f_\tau(\theta)}=\log\left(\frac{f_1(\theta)}{f_0(\theta)}\right),
\]
is independent of $\tau$. Therefore the following statement holds.

\begin{prop}{\sf The length of the logarithmic path connecting two power spectral densities $f_0$ and $f_1$, is given by
\begin{eqnarray}\label{eq:onemore}
\ell_{\rm logpath}:=\sqrt{
\intpi\left(\log\frac{f_1(\theta)}{f_0(\theta)}\right)^2d\theta- 
\left(\intpi\log\frac{f_1(\theta)}{f_0(\theta)}d\theta\right)^2
}
\end{eqnarray}
}
\end{prop}
\vspace*{.05in}

\begin{proof} The proof follows in the arguments leading to the proposition.\end{proof}

\section{Degradation of the smoothing error variance}\label{sec:III}

In a way completely analogous to the previous sections we may consider the increase in the variance of the smoothing error when a wrong choice between two alternatives is used to identify a candidate smoothing filter.

Thus, we begin with two density functions $f_1,f_2$ and assume that $f_1^{-1},f_2^{-1}\in L_1[-\pi,\pi)$. Accordingly we test the optimal smoothing filter based on $f_2$ against $f_1$.
As explained in Section \ref{section:preliminaries}, the $f_2$-optimal smoothing filter gives rise to an error
$u_0 -\hat{u}_{0|{\rm past \;\&\; future}}$ corresponding, via the Kolmogorov mapping, to
$h_{f_2}f_2(\theta)^{-1}$.
Hence, the variance of the smoothing error divided by the $f_1$-optimal variance is
\begin{eqnarray*}
{\rho}_{\rm smooth}(f_1,f_2) &:=&\left(\intpi | h_{f_2}f_2(\theta)^{-1}|^2 f_1(\theta)d\theta\right)\frac{1}{h_{f_1}}\\
&&\hspace*{-1cm}=\left(\intpi f_2(\theta)^{-2} f_1(\theta)d\theta\right)\frac{h_{f_2}^2}{h_{f_1}}\\
&&\hspace*{-1cm}=\left(\intpi \frac{f_1(\theta)}{f_2(\theta)^2}d\theta\right)\frac{\left(\intpi f_1(\theta)^{-1}d\theta\right)}
{\left(\intpi f_2(\theta)^{-1}d\theta\right)^2}
\end{eqnarray*}
Interestingly, this can be rewritten as follows
\begin{eqnarray}\nonumber
{\rho}_{\rm smooth}(f_1,f_2) &=&
\left(\intpi \left(\frac{f_1(\theta)}{f_2(\theta)}\right)^2 f_1(\theta)^{-1}d\theta\right)
\\&&\times \frac{\left(\intpi f(\theta)^{-1}d\theta\right)}
{\left(\intpi \left(\frac{f_1(\theta)}{f_2(\theta)}\right) f_1(\theta)^{-1}d\theta\right)^2}\nonumber\\
&=& \left(\intpi \left(\frac{f_1(\theta)}{f_2(\theta)}\right)^2 d\phi_1(\theta)\right)\nonumber
\\&&\times \frac{1}{\left(\intpi \left(\frac{f_1(\theta)}{f_2(\theta)}\right) d\phi_1(\theta)\right)^2}\nonumber\\
&=& \left(\frac{\sqrt{\intpi \left(\frac{f_1(\theta)}{f_2(\theta)}\right)^2 d\phi_1(\theta)}}
{\intpi \left(\frac{f_1(\theta)}{f_2(\theta)}\right) d\phi_1(\theta)}\right)^2
\label{eq:avsms}
\end{eqnarray}
where
\[
d\phi_1(\theta):=\frac{f_1(\theta)^{-1}d\theta}{\intpi f_1(\theta)^{-1}d\theta}
\]
is a normalized measure with variation one. Expression (\ref{eq:avsms}) shows the degradation as the square of the ratio
of the mean-square of the fraction $f_1/f_2$ over its arithmetic mean. These two means, mean-square and arithmetic, are weighted by $d\phi_1$ which is of course dependent on one of the two arguments. However, the expression is homogeneous and does not depend on scaling of either of the two arguments $f_1$ or $f_2$.

Accordingly, we may define as a distance measure
\begin{eqnarray}\nonumber
{\delta}_{\rm smooth}(f_1,f_2) &=&
\log(\rho_{\rm smooth}(f_1,f_2))\\
&=& \log\left({\intpi \left(\frac{f_1(\theta)}{f_2(\theta)}\right)^2 d\phi_1(\theta)}\right)
%\nonumber \\&&\hspace*{-1cm}
-\log\left(\left(\intpi \left(\frac{f_1(\theta)}{f_2(\theta)}\right) d\phi_1(\theta)\right)^2\right)
\label{eq:avsmslogs}
\end{eqnarray}
The presence of a data-dependent integration measure may by compared to the (normalized) Kullback-Leibler divergence in (\ref{eq:KL}).

\section{Reappraisal and generalizations}

The expressions derived in the previous sections suggest that generalized means of the ``likelihood''-like ratio $\Lambda:=f_1/f_2$ and their logarithms may be used as distance measures between ``shapes'' of density functions $f_1$ and $f_2$.  More specifically, we know that for any positive function $\Lambda$,
\[
M_r(\Lambda)\leq M_s(\Lambda) \mbox{ for any }-\infty\leq r<s\leq \infty,
\]
where $M_r(\Lambda)$ denotes the $r$-th generalized mean
\[
M_r(\Lambda):=\left(\intpi \Lambda(\theta)^rd\theta\right)^{1/r}.
\]
Then
\[ \delta_{r,s}(\Lambda):=\log(M_r(\Lambda))-\log(M_s(\Lambda)) \geq 0
\]
with a value which depends on how ``far'' $\Lambda$ is from being constant. Hence, we may use $\delta_{r,s}(f_1/f_2)$ to quantify the distance between the ``shapes'' of $f_1$ and $f_2$, and
since
\[
M_0(\Lambda):=\lim_{r\to 0}M_r(\Lambda)=e^{\intpi \log(\Lambda(\theta))d\theta}
\]
is the geometric mean of $\Lambda$ (see e.g., \cite[page 23]{BB}), both $\delta,\, \delta_{a/g}$ that we encountered earlier are special cases of the above. Weighted versions of weighted means may also be used for the same purpose, as suggested in Section \ref{sec:III}. Lengths of geodesics as suggested in Section \ref{sec:riemannian} present another possibility. Indeed, a ``zoo'' of possible options emerges. Assessing practical and theoretical merits of each is the subject of a future project.

%We conclude by noting that the conceptual framework of comparing degradation of performance with regard to specific tasks when the wrong choice between alternatives is used, extends naturally to other contexts. For instance, a similar rationale may be used in the context of positive-definite Toeplitz matrices and prediction based on a finite window of observations.


\begin{thebibliography}{99}
%\bibitem{dickinson} B.W. Dickinson, ``Two dimensional Markov spectrum estimates need not exist,''
%{\em IEEE Trans.\ on Information Theory}, {\bf IT-26}: 120-121, 1980.
\bibitem{alisilvey} S.M. Ali and S.D. Silvey, ``A general class of coefficients of divergence of one distribution from another,'' {\em J. Royal Stat.\ Soc.}, {\bf 28}: 131-142, 1966.

\bibitem{bhattacharyya} A. Bhattacharyya, ``On a measure of divergence between two statistical populations defined by their probability distributions,'' {\em Bull.\ Calcutta Math.\ Soc.}, {\bf 35}: 99-109, 1943.

\bibitem{BB} E.F. Beckenbach and R. Bellman, {\bf Inequalities}, Springer-Verlag, Berlin-Heidelberg, 198 pages,1965.

\bibitem{Bregman} L.M. Bregman, ``The relaxation method of finding the common point of  convex sets and its application to the solution of problems in convex programming,'' {\em USSR Comput.\ Math.\ and Math.\ Phys.} vol. 7, pp. 200-217, 1967.

%\bibitem{burg} J. Burg, {\bf Maximum entropy spectral analysis}, Ph.D. dissertation, Stanford University, 1975.

\bibitem{csiszar2} I. Csisz\'{a}r, ``Why least squares and maximum entropy? An axiomatic approach to inference for linear inverse problems,'' {\em The Annals of Probability,} {\bf 19(4)}: 2032-2066, 1991.

%\bibitem{dgk} P. Delsarte, Y. Genin, and Y. Kamp,
%``Half-Plane Toeplitz Systems,'' {\em IEEE Transactions on Information Theory,} {\bf IT-26}: 465-474, 1980.

%\bibitem{donoghue} W.F. Donoghue, Jr., {\bf Monotone Matrix functions and Analytic Continuation}, Springer-Verlag, 1974.

%\bibitem{sesha2} T.T. Georgiou, {\sf Spectral Estimation via Selective Harmonic
%Amplification,} {\em IEEE Trans.\ on Automatic Contr.},
%January 2001, {\bf 46(1)}: 29-42.

%\bibitem{acmatrix1} T.T. Georgiou,
%``The structure of state covariances and its relation to
%the power spectrum of the input,''
%{\em IEEE Trans.\ on Automatic Control}, {\bf 47(7):} 1056-1066, July 2002.

%\bibitem{acmatrix2}
%T.T. Georgiou, ``Spectral analysis based on the state covariance:
%the maximum entropy spectrum and linear fractional parameterization,''
%{\em IEEE Trans. on Automatic Control}, {\bf 47(11):}
%1811-1823, November 2002.

\bibitem{BGL2} C. Byrnes, T.T. Georgiou, and A. Lindquist,
``A new approach to spectral estimation: A tunable high-resolution spectral estimator,''
{\em IEEE Trans.\ on Signal Processing}, {\bf 48(11):} 3189-3206, November 2000.

\bibitem{ac_may2004} T.T. Georgiou,
``Solution of the general moment problem  via a one-parameter imbedding,''
{\em IEEE Trans.\ on Automatic Control}, {\bf 50(6):} 811-826, June 2005.

\bibitem{IT}
T.T. Georgiou,
``Relative Entropy and the multi-variable multi-dimensional Moment Problem,''
{\em IEEE Trans. on Information Theory}, {\bf 52(3):} 1052 - 1066,  March 2006. 

%\bibitem{CAR} T.T. Georgiou,
%``The Carath\'{e}odory-Fej\'{e}r-Pisarenko decomposition
%and its multivariable counterpart,'' preprint, 29 pages: arXiv:math.OC/0509225v1; http://arxiv.org/abs/math/0509225/

\bibitem{ansatz} T.T. Georgiou,
``The maximum entropy ansatz in the absence of a time arrow: fractional-pole models,''
preprint, 18 pages: http://arxiv.org/abs/math/0601648/
%\bibitem{car}
%T.T. Georgiou,
%``Decomposition of Toeplitz matrices via convex optimization,'' {\em IEEE Signal Processing Letters,} submitted December 2005.

%\bibitem{Geronimus} Ya. L. Geronimus,
%{\bf Orthogonal Polynomials}, English translation from Russian by Consultants Bureau,
%New York, 570 pages, 1961.

\bibitem{GrenanderSzego} U. Grenander and G. Szeg\"{o}, {\bf Toeplitz Forms and their
Applications}, Chelsea, 1958.

\bibitem{Haykin} S. Haykin, {\bf Nonlinear Methods of Spectral Analysis,}
Springer-Verlag, New York, 247 pages, 1979.

\bibitem{GL}
T.T. Georgiou and A. Lindquist, ``Kullback-Leibler approximation of spectral
density functions,'' {\em IEEE Trans. on Information Theory}, {\bf 49(11)},
November 2003.

\bibitem{Gray} R. Gray, A. Buzo, A. Gray, and Matsuyama, ``Distortion measures for speech processing,'' {\em IEEE Trans.\ on Acoustics, Speech, and Signal Proc.}, vol.\ 28, no.\ 4, Aug.\ 1980.

%\bibitem{Haykin_array} S. Haykin, J.H. Justice, N.L. Owsley, J.L. Yen, and A.C. Kak, {\bf Array Signal Processing}, Prentice-Hall, 1985.

%\bibitem{halliwell} J.J. Halliwell, J. P\'{e}rez-Mercader, and W.H. Zurek, Physical origins of time asymetry, Cambridge Univ.\ Press, 1994.

%\bibitem{jaynes} E.T. Jaynes, ``On the rationale of maximum entropy methods,'' {\em Proc.\ IEEE},
%{\bf 70}: 939-952, 1982.

%\bibitem{junk} M. Junk, ``Maximum entropy for reduced moment problems,'' {\em Mathematical Models and Methods in Applied Sciences}, {\bf 10(7)}: 1001-1025, 2000.

%\bibitem{kittel} C. Kittel, Thermal Physics, Wiley, 1969.

%\bibitem{KreinNudelman} M.G. Krein and A.A. Nudel'man,
%{\bf The Markov Moment Problem and Extremal Problems,}
%American Mathematical Society, Providence, RI, 417 pages, 1977.

%\bibitem{LevineTribus} R.D. Levine and M. Tribus (editors), {\bf The Maximum Entropy Formalism}, MIT Press, Cambridge, 1979.

\bibitem{masani} P. Masani,
{\sf Recent trends in multivariate prediction theory,} in {\bf Multivariate Analysis},
P.R. Krishnaiah, Ed., Academic Press, pp.\ 351-382, 1966.

\bibitem{Rudin} W. Rudin, {\bf Real and Complex Analysis}, 3rd edition,
McGraw Hill, 1987.

%\bibitem{rudin} W. Rudin, ``The extension problem for positive-definite functions,'' {\em Illinois J. Math.}, {\bf 7}:532-539, 1963.

%\bibitem{vantrees} H.L. Van Trees, {\bf Optimum Array Processing: part IV of Detection, Estimation and Modulation Theory}, Wiley-Interscience, 2002.

\bibitem{StoicaMoses} P. Stoica and R. Moses, {\bf Introduction to Spectral Analysis},
Prentice Hall, 2005.

\bibitem{Szego} G. Szeg\"{o}, ``\"{U}ber die randwerten eiher analytischen functionen,'' {\em Math.\ Ann.}, {\bf 84:} 232-244, 1921.

\bibitem{Varadhan} S.R.S. Varadhan, {\bf Probability Theory}, AMS, 2000.

\end{thebibliography}
\end{document}